\documentclass[10pt]{article}
\font\smallit=cmti10
\font\smalltt=cmtt10

\usepackage{amssymb,latexsym,amsmath,epsfig,amsthm} 
	\usepackage{hyperref}   
	\hypersetup{colorlinks=true, linkcolor=blue, 
		filecolor=magenta, urlcolor=red,}
\makeatletter

\renewcommand\section{\@startsection {section}{1}{\z@}
{-30pt \@plus -1ex \@minus -.2ex}
{2.3ex \@plus.2ex}
{\normalfont\normalsize\bfseries\boldmath}}

\renewcommand\subsection{\@startsection{subsection}{2}{\z@}
{-3.25ex\@plus -1ex \@minus -.2ex}
{1.5ex \@plus .2ex}
{\normalfont\normalsize\bfseries\boldmath}}

\renewcommand{\@seccntformat}[1]{\csname the#1\endcsname. }

\makeatother

\newtheorem{theorem}{Theorem}

\theoremstyle{definition}

\begin{document}

\begin{center}
\uppercase{\bf On the Brousseau sums} $\sum_{i=1}^n i^p F_i$
\vskip 20pt
{\bf Gregory Dresden}\\
{\smallit Department of Mathematics, Washington and Lee  University, Lexington, VA, USA}\\
{\tt dresdeng@wlu.edu}\\ 
\end{center}
\vskip 20pt

\centerline{\bf Abstract}

\noindent
We start with new convolution formulas for $F_n - n^p$ involving only the binomial coefficients. We then use those to find direct formulas for the sums $\sum_{i=1}^n i^p F_{n-i}$ and $\sum_{i=1}^n i^p F_i$, and we show how our formulas connect to  work in earlier papers by Ledin, Brousseau, Zeitlin, Adegoke, Shannon and Ollerton, and Kinlaw, Morris and Thiagarajan.

\pagestyle{myheadings}
\markright{\smalltt INTEGERS: 22 (2022)\hfill}
\thispagestyle{empty}
\baselineskip=12.875pt
\vskip 30pt 

\section{Introduction}
Our story begins in 1963 when Brother
Alfred Brousseau published the following problem 
\cite{Brousseau1} in what was just the second issue of the newly-founded {\em Fibonacci Quarterly}:
\begin{equation*}
\mbox{Sum: \ \ \ } \sum_{i=1}^n i^3 F_i.
\end{equation*}
Brousseau is using the common definition of the Fibonacci numbers as  $F_0 = 0$, $F_1 = 1$, and $F_n = F_{n-1} + F_{n-2}$ for $n>1$, and we will do the same.
A solution by Erbacher and Fuchs \cite{ErbacherFuchs} appeared the following year. They showed that
\begin{equation}\label{solution1}
\sum_{i=1}^n i^3 F_i \ = \ 
    (n^3 + 6n - 12)F_{n+2} + 
    (-3n^2 + 9n - 19)F_{n+3} + 50,
\end{equation}
and this is both completely true and completely unsatisfying: Why 50? Why  is there a cubic with $F_{n+2}$ but not 
 with $F_{n+3}$? What happens if we try to sum $i^4F_i$?

A few years later, we find follow-up articles by 
Brousseau \cite{Brousseau2}, Ledin \cite{Ledin}, and Zeitlin \cite{Zeitlin}, all of whom  place the original problem  in the proper context of finding 
solutions to
$\sum_{i=1}^n i^p F_i$ for arbitrary powers $p$.
Ledin re-writes the solution 
in Equation (\ref{solution1}) as 
\[
\sum_{i=1}^n i^3 F_i \ = \  
    (n^3 -3n^2 + 15n -31)F_{n} +
    (n^3-6n^2+24n-50)F_{n+1} + 50,
\]
and he notes that the coefficients $\{1,-3,15,-31\}$ and 
$\{1, -6, 24, -50\}$ seen above, and similar lists of coefficients for solutions to 
$\sum_{i=1}^n i^p F_i$, 
can be written as products of:
binomial coefficients, and a sequence of numbers ``...  the law of formation of which is yet to be determined". Ledin's work on finding solutions to $\sum i^p F_i$ involves 
a complicated set of four separate integrals.
Zeitlin, writing that same year, is able to determine that the ``law of formation" involves factorials, Fibonacci numbers, and Stirling numbers of the second kind.

Jumping ahead to the 2020's, we see a number of authors shedding some new light onto this old problem. 
Kinlaw, Morris and Thiagarajan \cite{KMT} show that
the coefficients $\{1,-3,15,-31\}$ and 
$\{1, -6, 24, -50\}$ can be written in terms of Eulerian numbers, and they develop a number of nice properties for them.  
Adegoke \cite{Adegoke} also defines them in terms of Eulerian numbers, and sets up some lovely recursive definitions using derivatives. 
Shannon and Ollerton \cite{ShannonOllerton} conjecture a fascinating formula relating these 
to the Bernoulli numbers.
We will return to all of these authors and their work in a moment.

And as for us? We start with 
the nice little identity
\begin{equation}\label{nicelittle}
F_n - n^3 \ = \ \sum_{i=1}^{n} \left(i^3 - 2\left( \binom{3}{2}i^2+\binom{3}{0}i^0\right)\right) \cdot F_{n-i}.
\end{equation}
We use that formula (along with similar ones for $F_n - n^2$ and $F_n - n$) to arrive at this convolution formula,
\begin{equation*}
\sum_{i=1}^{n} i^3 F_{n-i}  \ = \  
31 F_n + 50 F_{n+1} - (n^3 + 6n^2 + 24n + 50),
\end{equation*}
and the attentive reader will recognize those coefficients
$\{1, 6, 42, 50\}$ from earlier. 
From here, we can produce Ledin's solution to the original problem of $\sum_{i=1}^{n} i^3 F_{i}$, but this time with a clear and direct understanding of the coefficients. We give a recursive formula for these coefficients using only binomial coefficients and without any mention of either Stirling numbers or the Eulerian numbers. Likewise, we do not need to call upon the complicated integration techniques of Ledin or the clever derivative methods of Adegoke. 
Of course, our work applies not just to 
$\sum_{i=1}^{n} i^3 F_{i}$
but in fact to 
$\sum_{i=1}^{n} i^p F_{i}$ for all non-negative powers $p$, and along the way we establish some lovely and novel convolution formulas for $F_n - n^p$.

\section{Differences of Fibonacci numbers and powers}\label{section2}

So, with all this in mind, let us begin with our first set of formulas. We have mentioned that 
\[
F_n - n^3 = \sum_{i=1}^{n} \left(i^3 - 2\left( \binom{3}{2}i^2+\binom{3}{0}i^0\right)\right) \cdot F_{n-i},
\]
and now we put this in the proper context: 
	\begin{equation*} 
		\begin{aligned}
F_{n} - n\,\, \ &= \ \sum_{i=1}^{n}  (i-2(1)) \cdot F_{n-i}, \\
F_{n} - n^2 \ &= \ \sum_{i=1}^{n} (i^2 -2(2i)) \cdot F_{n-i}, \\
F_{n} - n^3 \ &= \ \sum_{i=1}^{n} (i^3 - 2(3i^2+1)) \cdot F_{n-i}, \\
F_{n} - n^4 \ &= \ \sum_{i=1}^{n} (i^4 - 2(4i^3 + 4i )) \cdot F_{n-i}. 
		\end{aligned}
	\end{equation*}
We note that the numbers on the right of each equation are indeed the binomial coefficients. The first
two equations appear on the OEIS 
\cite{OEIS} at
\href{https://oeis.org/A065220}{A065220} and 
\href{https://oeis.org/A014283}{A014283},
respectively. 
The general formula is given here. 

\begin{theorem}\label{Fn.minus.np}
For $p \geq 1$, we have 
\begin{equation}
F_n - n^{p} \ = \ \sum_{i=1}^{n} \left( i^{p} -   2\sum_{j=0}^{p/2} \binom{p}{2j+1} i^{p-2j-1} \right) \cdot F_{n-i}. \label{theorem1.1}
\end{equation}
\end{theorem}

Our proof is elementary, using only the properties of binomial coefficients and the recursion formula for the Fibonacci numbers.

\begin{proof}
We proceed by induction on $n$. For $n=1$, the left-hand side is $F_1 - 1^{p}$ which is zero, 
and the right-hand side is a single term multiplied by $F_{0}$ which is, again, zero.

For $n=2$ it takes a bit of work. On the left, we have $F_2 - 2^{p}$ which is $1-2^{p}$. 
On the right, we have 
\[
\sum_{i=1}^{2} \left(i^{p} \ - \ 2\sum_{j=0}^{p/2} \binom{p}{2j+1} i^{p-2j-1} \right) \cdot F_{2-i},
\]
and since $F_0 = 0$, the only non-zero part of this sum is at $i=1$, giving us 
\[
\left(1 \ - \ 2\sum_{j=0}^{p/2} \binom{p}{2j+1} \right) \cdot F_{2-1}.
\]
The sum on the inside is made up of every other term in the $p$th row of Pascal's
Triangle, and it is well-known that while the sum of every term is $2^{p}$, the sum 
of every other term is simply $2^{p-1}$. This gives us 
\[
\Bigg( 1 \ - \ 2 \cdot 2^{p-1}\Bigg) \cdot F_{2-1},
\]
which simplifies to $1-2^{p}$, as desired. 

Now, we do the induction step. We fix $n > 1$, and we assume that the formula holds for $n-1$ and $n$, and we attempt to prove that it holds for $n+1$. We will need the following identity, which is easily verified:
\begin{equation}\label{nnn}
(n-1)^{p} + n^{p} - (n+1)^{p} \ = \ n^{p} - 2 \sum_{j=0}^{p/2} \binom{p}{2j+1}n^{p-2j-1}.
\end{equation}
Next, we write out  two versions of Equation (\ref{theorem1.1}), at $n-1$ and $n$, respectively:
\begin{align}
F_{n-1} - {(n-1)}^{p}\,\,\,\,\, \ &= \ \sum_{i=1}^{n-1} \left(i^{p} \ - \ 2\sum_{j=0}^{p/2} \binom{p}{2j+1} i^{p-2j-1} \right) \cdot F_{n-1-i}, \label{n1} \\
F_n - n^{p}\,\,\,\,\, \ &= \ \sum_{i=1}^{n} \left(i^{p} \ - \ 2\sum_{j=0}^{p/2} \binom{p}{2j+1} i^{p-2j-1} \right) \cdot F_{n-i}. \label{n2}
\end{align}
Now, we add together Equations (\ref{nnn}), (\ref{n1}), and (\ref{n2}). On the left, we have $F_{n-1} + F_{n} - (n+1)^{p}$, which simplifies to 
$F_{n+1} - (n+1)^{p}$, as desired. On the right, the two sums in Equations (\ref{n1}) and (\ref{n2}) differ only in the upper limit of the outer sum, and in the Fibonacci number on the right of each equation. Hence, when we add together these two sums, we get 
\begin{equation}\label{nn3}
\ \ \ \ \ \ \ \ \ \ \ \ \ \ \ \ \  \ \sum_{i=1}^{n-1} \left(i^{p} \ - \ 2\sum_{j=0}^{p/2} \binom{p}{2j+1} i^{p-2j-1} \right) \cdot F_{n+1-i},
\end{equation}
and this is because the $i=n$ term in Equation (\ref{n2}) is zero thanks to $F_{n-i}$ and so can be dropped. We now add 
the right-hand side of Equation (\ref{nnn}) to Equation (\ref{nn3}), recognizing that the right-hand side
of Equation (\ref{nnn}) is simply the missing $i=n$ term from the sum in Equation (\ref{nn3}). This gives us 
\begin{equation}\label{nn4}
\ \ \ \ \ \ \ \ \ \ \ \ \ \ \ \ \  \ \sum_{i=1}^{n} \left(i^{p} \ - \ 2\sum_{j=0}^{p/2} \binom{p}{2j+1} i^{p-2j-1} \right) \cdot F_{n+1-i}
\end{equation}
as the sum of the right-hand sides of Equations (\ref{nnn}), (\ref{n1}), and (\ref{n2}). Finally, we recognize that we can replace the $n$ in the upper limit of the sum in Equation (\ref{nn4}) with $n+1$, since the $i=n+1$ term is zero thanks to $F_{n+1-i}$. Gathering together all of our work, we see that the sum of Equations (\ref{nnn}), (\ref{n1}), and (\ref{n2}), looking at the right-hand sides and left-hand sides separately, does indeed give us
\[
F_{n+1} - {(n+1)}^{p}\,\,\,\,\, \ = \ \sum_{i=1}^{n+1} \left(i^{p} 
\ - \ 2\sum_{j=0}^{p/2} \binom{p}{2j+1} i^{p-2j-1} \right) \cdot F_{n+1-i},
\]
and this completes the induction step.
\end{proof}


\section{Convolutions}

Let us now turn things around. Recall that in Theorem \ref{Fn.minus.np} 
the formulas in Equation (\ref{theorem1.1})
had nice ``solutions" (like $F_n -n^2$ and $F_n - n^3$) but complicated ``summands". This next group has nice summands but complicated solutions that, at first glance, do not seem to follow any discernible pattern. Notice that we are now starting our sums at $i=0$ instead of at $i=1$:
	\begin{equation*} 
		\begin{aligned}
\sum_{i=0}^{n}   1 \cdot F_{n-i}  \  &=  \   1 \,F_n +  1 \, F_{n+1}  - 1, \\
\sum_{i=0}^{n}   i \cdot F_{n-i}  \  &=  \   1 \, F_n +  2 \, F_{n+1}  - (n+2), \\
\sum_{i=0}^{n} i^2 \cdot F_{n-i}  \  &=  \   5 \, F_n +  8 \, F_{n+1}  - (n^2+4n+8), \\
\sum_{i=0}^{n} i^3 \cdot F_{n-i}  \  &=  \  31 \, F_n + 50 \, F_{n+1}  - (n^3+6n^2+24n + 50), \\
\sum_{i=0}^{n} i^4 \cdot F_{n-i}  \  &=  \ 257 \, F_n +416 \, F_{n+1}  - 
		(n^4 + 8n^3 + 48n^2 + 200n + 416).
	\end{aligned}
	\end{equation*}
Naturally, there is indeed a pattern to these solutions. To understand the pattern, we must first define the following two sequences of numbers. For our first sequence, we define $A_0 = 1$ and then we define $A_p$ for $p>0$ as
	\begin{equation} \label{AGroup}
A_{p}   \ = \ (-1)^p + 2\sum_{j=0}^{p/2} \binom{p}{2j+1} A_{p-2j-1}. 
	\end{equation}
This gives us the numbers $1,1,5,31,257,\dots$ which are the coefficients of $F_n$ in each line of the five equations given above.  This sequence is \href{https://oeis.org/A000556}{A000556} on the OEIS. 

For our second sequence, we define $B_0 = 1$ and then we define $B_p$ for $p>0$ as
	\begin{equation} \label{BGroup}
B_{p}   \ = \  2\sum_{j=0}^{p/2} \binom{p}{2j+1} B_{p-2j-1}. 
	\end{equation}
This will produce the numbers $1, 2, 8, 50, 416, \dots$,
 which appear twice in each line of the five equations given above, first as the 
 coefficients of $F_{n+1}$  and again as the constant terms. This  sequence is  \href{https://oeis.org/A000557}{A000557} on the OEIS. (We will reveal later that these two sequences $A_p$ and $B_p$ are equal to the sequences $M_{1,p}$ and $M_{2,p}$ 
 as used by Zeitlin, Ledin, and others.)
 
 To further understand the patterns in 
 the set of five equations given above,
  we also need to recognize that each sum on the left
  of 
 those equations is actually a {\em convolution} of 
 the Fibonacci numbers and the powers of $i$. To be precise, we define this convolution, $\mathcal{C}_n^{(p)}$, as follows:
 \begin{equation}\label{convolutiondefinition}
  \mathcal{C}_n^{(p)} = 
  	\left\{ \begin{array}{rl}
  	            \displaystyle\sum_{i=0}^{n} \, 1 \cdot F_{n-i}  \, & \mbox{if $p=0$},\\[4.0ex]
				\displaystyle\sum_{i=0}^{n} \, i^p \cdot F_{n-i}  & \mbox{if $p>0$}.
			\end{array} \right.
 \end{equation}

With this in mind, here is our next theorem. 
\begin{theorem}\label{iFi.theorem}
For $\mathcal{C}_n^{(p)}$ the convolution as defined 
in Equation (\ref{convolutiondefinition}),
we have 
\begin{equation}\label{iFi.equation}
\mathcal{C}_n^{(p)} \ = \ A_p F_n + B_p F_{n+1} - 
\sum_{k=0}^p \binom{p}{k} B_k n^{p-k}.
\end{equation}
\end{theorem}

This is, quite simply, a stunningly beautiful theorem that gives order to the 
chaos of equations and constants in 
the five equations at the beginning of this section.
Our proof, once again, uses only elementary methods, although the details are rather technical. 

\begin{proof}
We will use induction on $p$ to show that 
Equation (\ref{iFi.equation}) holds for 
all $p \geq 0$, thus proving our theorem.
For the base case of $p=0$, it is easy to show that 
$\mathcal{C}_n^{(0)} = 1F_n + 1F_{n+1} - 1$, matching nicely with Theorem 
\ref{iFi.theorem} and with the first equation in  the group of five at the beginning of this section.

Before proceeding to the induction case of $p \geq 1$, we first need 
to re-write Equation (\ref{theorem1.1}) in Theorem \ref{Fn.minus.np}. This equation, which holds for $p \geq 1$, can be expanded to give us
\begin{equation}\label{sum1.no}
 F_n - n^p 
 \ = \ \sum_{i=1}^{n}  i^{p} F_{n-i}  \ - \   2\sum_{j=0}^{p/2} \binom{p}{2j+1} \sum_{i=1}^n i^{p-2j-1} 
    F_{n-i}.
\end{equation}
This is almost, but not quite, in the format we need, as these sums start at $i=1$ but to write these as convolutions we want our sums to start at $i=0$. If we do so, we will obtain the following formula:
\begin{equation}\label{sum1.yes}
 (-1)^p F_n - n^p 
 \ \ = \ \ \mathcal{C}_n^{(p)}   \ - \   2\sum_{j=0}^{p/2} \binom{p}{2j+1} \mathcal{C}_n^{(p-2j-1)}.
\end{equation}
Let us explain how this works. If $p\geq 1$ is even, 
then we can definitely add the $i=0$ terms of $0^pF_{n-0} = 0$ to the first sum
in Equation (\ref{sum1.no}),
and $0^{p-2j-1}F_{n-0} = 0$ to the last  sum in Equation  (\ref{sum1.no}), converting those sums into convolutions and giving us the right-hand side of 
Equation (\ref{sum1.yes}). Likewise, the left-hand sides of 
Equations (\ref{sum1.no}) and  (\ref{sum1.yes}) 
also match for 
$p \geq 1$ even, giving us that Equation (\ref{sum1.yes}) does indeed
hold for $p$ even. 

For $p$ odd, it is a bit harder. We can still add 
$0^pF_{n-0} = 0$ to the first sum
in Equation  (\ref{sum1.no}), converting it to the convolution
$\mathcal{C}_n^{(p)}$. To convert the last sum in
Equation (\ref{sum1.no}) to a convolution, we 
note that at $j=\lfloor p/2 \rfloor$, the last sum is 
actually $\sum_{i=1}^n i^0 F_{n-i}$ and so the ``missing" $i=0$ term
 is $1\,F_n$.
If we add and subtract that $1 \, F_n$ inside the sum on the right of Equation (\ref{sum1.no}) for $j=\lfloor p/2 \rfloor$,
converting all those sums to convolutions, we obtain
\begin{equation}\label{sum1.yes2}
 F_n - n^p 
 \ \ = \ \ \mathcal{C}_n^{(p)}   \ - \   2\sum_{j=0}^{p/2} \binom{p}{2j+1}  \mathcal{C}_n^{(p-2j-1)} \ + 2 \binom{p}{2 \lfloor p/2 \rfloor + 1}F_n \qquad \mbox{for $p$ odd}.
\end{equation}
That last term is simply $2F_n$, and so if we subtract it from both sides of Equation (\ref{sum1.yes2}) we will get 
\begin{equation}\label{sum1.yes3}
 - F_n - n^p 
 \ \ = \ \ \mathcal{C}_n^{(p)}   \ - \   2\sum_{j=0}^{p/2} \binom{p}{2j+1}  \mathcal{C}_n^{(p-2j-1)}  \qquad \mbox{for $p$ odd}.
\end{equation}
Of course, this just tells us that  Equation  (\ref{sum1.yes}) does indeed hold for $p$ odd.

We conclude that Equation (\ref{sum1.yes}) is true for 
all $p\geq 1$, and so we can now utilize Equation (\ref{sum1.yes}) to 
establish Equation (\ref{iFi.equation}) in the statement
of Theorem \ref{iFi.theorem}. To do so, let us re-write
Equation (\ref{sum1.yes}) in terms of 
$\mathcal{C}_n^{(p)}$, giving us a nice recursive definition 
for our convolutions:
\begin{equation}\label{iFi.equation3}
\mathcal{C}_n^{(p)} \ = \ (-1)^pF_n - n^{p} + 2\sum_{j=0}^{p/2} \binom{p}{2j+1} \mathcal{C}_n^{(p-2j-1)}.
\end{equation}
We will use this recursive description of 
$\mathcal{C}_n^{(p)}$, along with our induction hypothesis that  Equation (\ref{iFi.equation}) holds for all values less than $p$, to show that the coefficient of $F_n$ in Equation (\ref{iFi.equation}) really is $A_p$, that the coefficient of $F_{n+1}$ really is $B_p$, 
and that the coefficients of the powers of $n$ in Equation (\ref{iFi.equation}) 
really are $-\binom{p}{k}B_k$.

First, the coefficient of $F_n$. From 
Equation (\ref{iFi.equation3}), we have a single $(-1)^p F_n$ out front, and then each $\mathcal{C}_n^{(p-2j-1)}$ in
the summand of Equation (\ref{iFi.equation3}) contributes $A_{p-2j-1}F_{n}$ thanks to our induction hypothesis from Equation (\ref{iFi.equation}), and these are each multiplied by $2\binom{p}{2j+1}$ in Equation (\ref{iFi.equation3}) and summed over $j$. In other words, the $F_n$ terms in 
Equation (\ref{iFi.equation3}) are
\[
(-1)^p F_n \ + \  2 \sum_{j=0}^{p/2} \binom{p}{2j+1}
A_{p-2j-1}F_{n},
\]
and if we factor out the $F_n$ we have
\[
\left((-1)^p  \ + \  2 \sum_{j=0}^{p/2} \binom{p}{2j+1}
A_{p-2j-1}\right) \cdot F_n.
\]
From our definition in Equation (\ref{AGroup}) we see that this last equation is simply 
\[
A_p \cdot F_n,
\]
as desired. 

Next, the coefficient of $F_{n+1}$. We look again at 
equation  (\ref{iFi.equation3}), and we see that each $\mathcal{C}_n^{(p-2j-1)}$ in
the summand of Equation (\ref{iFi.equation3}) contributes $B_{p-2j-1}F_{n+1}$ thanks to our induction
hypothesis in Equation (\ref{iFi.equation}), and these are each multiplied by $2\binom{p}{2j+1}$ in Equation (\ref{iFi.equation3}) and summed over $j$. In other words, the $F_{n+1}$ terms in Equation  (\ref{iFi.equation3}) are
\[
2 \sum_{j=0}^{p/2} \binom{p}{2j+1}
B_{p-2j-1}F_{n+1},
\]
and if we factor out the $F_{n+1}$ we have
\[
\left(2 \sum_{j=0}^{p/2} \binom{p}{2j+1}
B_{p-2j-1}\right) \cdot F_{n+1}.
\]
From our definition in Equation (\ref{BGroup}) we see that this is simply 
\[
B_p \cdot F_{n+1},
\]
as desired.

Finally, the coefficients of the powers of $n$. We look once more at
Equation (\ref{iFi.equation3}), and we notice that there is a single $-n^p$ out front, and then each 
$\mathcal{C}_n^{(p-2j-1)}$ in
the summand of Equation (\ref{iFi.equation3}) contributes 
$- \sum_{k=0}^{p-2j-1} \binom{{p-2j-1}}{k} B_k\, n^{p-2j-1-k}$
thanks again to our induction hypothesis with Equation (\ref{iFi.equation}), and these are each multiplied by $2\binom{p}{2j+1}$ in Equation (\ref{iFi.equation3}) and summed over $j$. In other words, the 
``polynomial" part of Equation (\ref{iFi.equation3}) is
\[
-n^p - 2 \sum_{j=0}^{p/2} \binom{p}{2j+1} 
\sum_{k=0}^{p-2j-1} \binom{{p-2j-1}}{k} B_k \, n^{p-2j-1-k},
\]
and moving the second sigma from the inside to the outside, we have 
\begin{equation}\label{doublesum1}
-n^p - 2 \sum_{j=0}^{p/2} 
\sum_{k=0}^{p-2j-1} 
\binom{p}{2j+1} \binom{{p-2j-1}}{k} B_k \, n^{p-2j-1-k}.
\end{equation}

Let us now look carefully at the ``domain" of this double sum. Table \ref{table:1} gives all the ordered pairs $(j,k)$ that have non-zero contributions to our double sum in 
Equation (\ref{doublesum1}). To avoid cluttering up the table, we are only writing the ``$n$" term. Note that the ``$*$" at position $(j,k) = 
(\lfloor p/2 \rfloor, 0)$ is zero if $p$ is even, and is $n^0$ if $p$ is odd. 

\begin{table}[h!]
\centering
\begin{tabular}{l| c c c c c} 
\footnotesize $k=p-1$ & $n^0$ & & & & \\
\footnotesize $k=p-2$ & $n^1$ & & & & \\
\footnotesize $k=p-3$ & $n^2$ & $n^0$ & & & \\
\footnotesize $k=p-4$ & $n^3$ & $n^1$ & & & \\
\footnotesize $k=p-5$ & $n^4$ & $n^2$ & $n^0$ & & \\
$\vdots$ & $\vdots$ & $\vdots$ & $\vdots$ & & \\
\footnotesize $k=4$ & $n^{p-5}$ & $\vdots$ & $\vdots$ & $\vdots$ &\\
\footnotesize $k=3$ & $n^{p-4}$ & $\vdots$ & $\vdots$ & $\vdots$ &  \\
\footnotesize $k=2$ & $n^{p-3}$ & $n^{p-5}$ & $\vdots$ & $\vdots$ & \\
\footnotesize $k=1$ & $n^{p-2}$ & $n^{p-4}$ & $\vdots$ &
$\vdots$ & \\[1.05ex]
\footnotesize $k=0$ & $n^{p-1}$ & $n^{p-3}$ & $n^{p-5}$& $\dots$ & $*$ \\
\hline
  &\footnotesize $j=0$ &\footnotesize $j=1$ &\footnotesize $j=2$ & $\dots$ & \footnotesize $j= \lfloor p/2 \rfloor $  
\end{tabular}
\caption{Ordered pairs, and powers of $n$, that appear in our double sum.}
\label{table:1}
\end{table}
    
We now execute the following change of variables: we keep $j$ as $j$ but let $K = 2j + 1 + k$, chosen precisely because with this substitution,
the expression $n^{p-2j-1-k}$ in Equation (\ref{doublesum1}) becomes
the much simpler $n^{p-K}$. 
Note also that the fixed values of $K$ correspond to steep downward diagonals in Table \ref{table:1}, such that when $K=1$ we have $(j,k)=(0,0)$ with $n^{p-1}$, and when $K=2$ we have 
$(j,k) = (0,1)$ with $n^{p-2}$, and when $K=3$ we have $(j,k) = (0,2)$ and $(1,0)$ both with $n^{p-3}$, and so on. 

From looking at the table, we see that $K$ runs from $1$ to $p$, and for a given $K$, we see that $j$ runs from
$0$ to $\lfloor K/2 \rfloor$. 
Since $K = 2j+1+k$, then when we re-write Equation (\ref{doublesum1}) in terms of $j$ and $K$, we replace 
the $k$ in Equation  (\ref{doublesum1}) with 
$k=K-2j-1$, giving us
\begin{equation}\label{doublesum2}
-n^p - 2 \sum_{K=1}^{p} 
\sum_{j=0}^{K/2} 
\binom{p}{2j+1} \binom{{p-2j-1}}{K-2j-1} B_{K-2j-1} \, n^{p-K}.
\end{equation}
We now have to play with those binomial coefficients in 
Equation (\ref{doublesum2}) just a bit, to get them in the right form. We use the identities 
$\binom{a}{b} = \binom{a}{a-b}$ and 
$\binom{a}{b} \binom{b}{c} = \binom{a}{c}\binom{a-c}{a-b}$  respectively, to obtain
\begin{equation*}
    \binom{p}{2j+1} \binom{{p-2j-1}}{K-2j-1} \ = \ 
    \binom{p}{p-2j-1} \binom{{p-2j-1}}{p-K} \ = \  
    \binom{p}{p-K} \binom{{K}}{2j+1},
\end{equation*}
and when we substitute this into Equation (\ref{doublesum2}) we obtain
\begin{equation}\label{doublesum3}
-n^p - 2 \sum_{K=1}^{p} 
\sum_{j=0}^{K/2} 
\binom{p}{p-K} \binom{K}{2j+1} B_{K-2j-1} \, n^{p-K}.
\end{equation}
We can pull out that first binomial coefficient, and bring in the $2$, to get
\begin{equation}\label{doublesum4}
-n^p - \sum_{K=1}^{p} 
\binom{p}{p-K}\left(2 \sum_{j=0}^{K/2} 
 \binom{K}{2j+1} B_{K-2j-1}\right) \, n^{p-K}.
\end{equation}
This inside sum is a perfect match for Equation (\ref{BGroup}), allowing us to replace it with $B_K$. We can also replace $\binom{p}{p-K}$ with $\binom{p}{K}$, and so 
Equation (\ref{doublesum4}) is now
\begin{equation}\label{doublesum5}
-n^p - \sum_{K=1}^{p} 
\binom{p}{K} \, B_K \, n^{p-K}.
\end{equation}
Since $B_0$ is defined as equal to $1$, we can bring in the $-n^p$ by starting the sum at $K=0$ instead of $K=1$. This gives us, finally, 
\begin{equation}\label{doublesum6}
 - \sum_{K=0}^{p} 
\binom{p}{K} \, B_K \, n^{p-K}
\end{equation}
as the sum of all the $n$ terms in Equation (\ref{iFi.equation3}), and this is indeed a perfect match for the expression in
Equation (\ref{iFi.equation}). This, along with our earlier work on $A_p F_n$ and $B_p F_{n+1}$, shows that Equation 
(\ref{iFi.equation}) really does hold for all $p \geq 0$,
thus concluding our proof.
\end{proof}


\section{The Brousseau sums}

Where do we go from here? Recall that our formulas in Theorem \ref{Fn.minus.np}
had nice ``solutions" (like $F_n -n^2$ and $F_n - n^3$) but complicated ``summands".
As for our formulas in Theorem \ref{iFi.theorem},
the solutions were rather complicated but the summands were simple convolutions of the Fibonacci numbers with the powers of integers. We can 
now bring everything together in this next group, which replaces the {\em convolutions} of Theorem \ref{iFi.theorem} with 
Brousseau's
{\em weighted sums}, as follows:
	\begin{equation*}
		\begin{aligned}
\sum_{i=0}^{n}   1 \cdot F_{i}  \  &=  \   (1)F_n \ + \   (1)  F_{n+1}  \ - \ 1, \\
\sum_{i=0}^{n}   i \cdot F_{i}  \  &=  \   (n-1)  F_n \ + \  (n-2)  F_{n+1} \ + \ 2, \\
\sum_{i=0}^{n} i^2 \cdot F_{i}  \  &=  \   (n^2-2n+5)  F_n \ + \  (n^2-4n+8)  F_{n+1} \  -\  8, \\
\sum_{i=0}^{n} i^3 \cdot F_{i}  \  &=  \  (n^3-3n^2+15n-31)  F_n \ + \ (n^3-6n^2+24n-50)  F_{n+1}  \ + \ 50. 
	\end{aligned}
	\end{equation*}

 If we define the
sum $\mathcal{S}_n^{(p)}$ of powers 
and Fibonacci numbers as 
\begin{equation}\label{sumdefinition}
  \mathcal{S}_n^{(p)} = 
  	\left\{ \begin{array}{rl}
		\displaystyle\sum_{i=0}^{n} \, 1 \cdot F_{i}  \, & \mbox{if $p=0$},\\[4.0ex]
				\displaystyle\sum_{i=0}^{n} \, i^p \cdot F_{i}  & \mbox{if $p>0$},
			\end{array} \right. 
 \end{equation}
 then we have the following remarkable theorem that bears a striking resemblance to Theorem \ref{iFi.theorem}. 
 
 \begin{theorem}\label{sum.theorem}
For $\mathcal{S}_n^{(p)}$ the weighted sums as defined 
in Equation (\ref{sumdefinition}), and with $A_k$ and $B_k$ as defined in Equations (\ref{AGroup}) and (\ref{BGroup}), 
we have 
\begin{equation}\label{sum.equation}
\mathcal{S}_n^{(p)} \ = \ 
\left(\sum_{k=0}^p \binom{p}{k} (-1)^k A_k n^{p-k}\right) F_n
\ + \ 
\left(\sum_{k=0}^p \binom{p}{k} (-1)^k B_k n^{p-k}\right) F_{n+1} 
\ - \ (-1)^p B_p.
\end{equation}
\end{theorem}
As we discussed earlier, Ledin, Zeitlin, and others all had similar versions of this theorem, but
none of them recognized that these coefficients
$A_k$ and $B_k$ (which they called $M_{1,k}$ and $M_{2,k}$) come from our simple recursion formulas in 
Equations (\ref{AGroup}) and (\ref{BGroup}),
using only the binomial coefficients.

\begin{proof} 
No need for induction in this proof;
thanks to Theorem \ref{iFi.theorem}, we can 
prove this directly if we apply some clever manipulations to our sums. We note that when $p=0$, the right-hand side of Equation (\ref{sum.equation}) is simply
$(A_0)F_n + (B_0)F_{n+1} - B_0$, which is 
equal to $F_n + F_{n+1} -1$, which is indeed equal to 
$\sum_{i=0}^{n} \, 1 \cdot F_{i}$ as desired.

For $p>0$,  we have
\begin{align*}
    \mathcal{S}_n^{(p)} \ &= \ \sum_{i=0}^{n} \, i^p \cdot F_{i},\\
\intertext{and replacing $i$ with $n-i$ we have}
    \mathcal{S}_n^{(p)} \ &= \ \sum_{i=0}^{n} \, (n-i)^p \cdot F_{n-i}
         \ = \ \sum_{i=0}^{n} \left( \sum_{k=0}^p \binom{p}{k} n^{p-k}(-1)^k i^k\right) \cdot F_{n-i}.
\end{align*}
We should note that when $i=0$ and $k=0$ the $i^k$ term on the inside
of the sum should be interpreted as being equal to $1$; this is because that $i^k$  comes from the $k$th term in the binomial expansion of $(n-i)^p$. 

We switch the order of summation, and simplify slightly, to obtain
\begin{equation*}
 \mathcal{S}_n^{(p)} \ = \ \sum_{k=0}^{p} \binom{p}{k} 
            (-1)^k  \left( \sum_{i=0}^n   i^k F_{n-i} \right) n^{p-k},
\end{equation*}
and with our understanding that the $i^k$ term is $1$ when $i=0$ and $k=0$ then this simplifies even further to
\begin{equation}\label{star1}
\mathcal{S}_n^{(p)} \ = \ \sum_{k=0}^{p} \binom{p}{k} 
            (-1)^k  \, \mathcal{C}_n^{(k)} \, n^{p-k}.
\end{equation}
From Equation (\ref{iFi.equation}) in Theorem \ref{iFi.theorem},
we can write $\mathcal{C}_n^{(k)}$ as
\begin{equation}\label{star2}
\mathcal{C}_n^{(k)} \ = \ A_k F_n + B_k F_{n+1} - 
\sum_{j=0}^k \binom{k}{j} B_j n^{k-j}.
\end{equation}
Let us now stitch together Equations (\ref{star1}) and (\ref{star2}), collecting the $F_n$ terms, the $F_{n+1}$ terms, and the
``constant" terms to establish Equation (\ref{sum.equation}) in the
statement of our theorem. 

First, we have the $F_n$ terms. From Equations (\ref{star1}) and (\ref{star2}), 
we see that the $F_n$ term is 
\[
\sum_{k=0}^p \binom{p}{k} (-1)^k \cdot A_k \, F_n  \cdot n^{p-k},
\]
which is a perfect match for the $F_n$ term in 
Equation (\ref{sum.equation}).
Likewise, 
the $F_{n+1}$ term is 
\[
\sum_{k=0}^p \binom{p}{k} (-1)^k \cdot B_k \, F_{n+1}  \cdot n^{p-k},
\]
which likewise matches the $F_{n+1}$ term in Equation (\ref{sum.equation}).
Finally, we look at the ``constant" terms which we hope will simplify to just 
$-(-1)^p B_p$.
From Equations  (\ref{star1}) and 
 (\ref{star2}), we factor out the negative to get the ``constant" expression
\begin{equation*}
-\sum_{k=0}^{p} \binom{p}{k} (-1)^k 
    \left( \sum_{j=0}^k \binom{k}{j} B_j n^{k-j} \right)
    n^{p-k}.
\end{equation*}
We carefully switch the order of summation, and simplify the inside, to obtain
\begin{equation*}
-\sum_{j=0}^{p} \sum_{k=j}^p  (-1)^k 
    \binom{p}{k}  \binom{k}{j} B_j n^{p-j},
\end{equation*}
and if we move around a few terms then this becomes
\begin{equation}\label{thirtysix}
-\sum_{j=0}^{p} B_j n^{p-j} \sum_{k=j}^p  (-1)^k 
    \binom{p}{k}  \binom{k}{j}.
\end{equation}
Next, we use a common binomial identity to write
$\binom{p}{k}  \binom{k}{j} = 
 \binom{p}{j}  \binom{p-j}{k-j}$, and if we move that 
 $\binom{p}{j}$ outside of the inner sum in 
Equation  (\ref{thirtysix}) then we obtain
\begin{equation*}
-\sum_{j=0}^{p} B_j n^{p-j}\binom{p}{j} \sum_{k=j}^p  (-1)^k 
      \binom{p-j}{k-j}.
\end{equation*}
We assign $K=k-j$, allowing us to re-write that inner sum in 
terms of $K$ instead of $k$, giving us 
\begin{equation}\label{thirtyeight}
-\sum_{j=0}^{p} B_j n^{p-j}\binom{p}{j} \sum_{K=0}^{p-j}  (-1)^{K+j} 
      \binom{p-j}{K}.
\end{equation}
Thanks to \cite[Identity 167]{BQ}, the inner sum is zero whenever $p-j>0$, so the entire expression in 
Equation (\ref{thirtyeight}) is zero
unless $j=p$, and when $j$ does equal $p$ then Equation (\ref{thirtyeight})
simplifies to 
\begin{equation}\label{thirtynine}
- B_p n^{p-p}\binom{p}{p} \sum_{K=0}^{0}  (-1)^{K+p} 
      \binom{p-p}{K},
\end{equation}
and this is just $-B_p (-1)^p$, as desired.
\end{proof}

\section{Connections to previous work}
To summarize our work, we defined our numbers $A_k$ and $B_k$ 
in Equations 
(\ref{AGroup}) and (\ref{BGroup}) in terms
of just the binomial coefficients, and we then showed that these numbers appear in the convolution formulas 
$\sum i^p F_{n-i}$ of Theorem \ref{iFi.theorem} 
which then led to their appearance in
the sum formulas 
$\sum i^p F_{i}$ of Theorem \ref{sum.theorem}.

A version of our convolution formula in Theorem \ref{iFi.theorem} can be discerned from Ledin's Theorem 4 if we replace his $n$ with $n-1$, although Ledin uses a third sequence of numbers that he calls $M_{3,k}$ that he defines in terms of $M_{1,k}$ and 
$M_{2,k}$ (which are our numbers $A_k$ and $B_k$ as mentioned earlier). We feel that this obscures the beauty of Theorem 
\ref{iFi.theorem}, and also prevents the reader from easily
discovering Equation (\ref{nicefact}) mentioned below. 

If we compare our Theorem \ref{sum.theorem} with Ledin's Equations (1), (3a), and (3b) in \cite{Ledin}, we see that
our $A_k$ and $B_k$ are identical to Ledin's 
$M_{1,k}$ and $M_{2,k}$, the exact numbers that Ledin was referring to when he said, ``the law of formation of which is yet to be determined". The difference between our work and Ledin's is that while we first defined $A_k$ and $B_k$ 
in terms of the recursion formulas 
in Equations (\ref{AGroup}) and (\ref{BGroup}) and then showed they
appear in the sum formulas $\sum i^p F_{i}$,
Ledin (and almost all of the subsequent authors) first defined what we call the
numbers $A_k$ and $B_k$ 
as the numbers that appear in the sum
formulas $\sum i^p F_i$, and then  found some nice properties for these numbers. In some cases, they worked with the polynomials 
in $n$ that appear in Equation (\ref{sum.equation}) rather than the individual coefficients $A_k$ and $B_k$, but the idea was the same. 

For the sake of completeness, we will list a few of these properties for $A_k$ and $B_k$, 
rewritten in terms of our notation.
First, we have this nice set from Ledin, which establishes a relationship between $A_k$ and $B_k$:
\begin{align}
    (-1)^k A_k \ &= \ \sum_{j=0}^k \binom{k}{j}B_j (-1)^j, \label{LedinA} \\
           B_k \ &= \ \sum_{j=0}^k \binom{k}{j}A_j. \label{LedinB}
\end{align}
This is a good place to mention that our Theorem \ref{iFi.theorem} gives us a nice ``inverse" of Ledin's  Equation (\ref{LedinB}); our new ``inverse" equation is
\begin{equation}\label{nicefact}
A_k \ = \ \sum_{j=0}^{k-1}
    \binom{k}{j} B_j.
\end{equation} 
To see this, simply take $n=1$ in Equation (\ref{iFi.equation}) 
from our Theorem \ref{iFi.theorem} to obtain
\[
0^p\cdot F_1 + 1^p\cdot F_0 \ = \ A_p F_1 + B_p F_2 - \sum_{k=0}^p \binom{p}{k} B_k, 
\]
and if we simplify this a bit, and replace $k$ with $j$ and $p$ with $k$, we get our desired equation,
\begin{equation} \label{DresdenA}
A_k = \sum_{j=0}^{k-1} \binom{k}{j} B_j.
\end{equation}
Note that if we carefully add Equations (\ref{LedinA}) and (\ref{DresdenA}), 
we can obtain a formula for $A_k$ that is strikingly similar to Equations (\ref{AGroup}) and (\ref{BGroup}); we leave the details to the reader.

Moving on, Ledin credits V.~E.~Hoggatt Jr.~for these lovely formulas:
\begin{align}
     A_k \ &= \ \sum_{j=0}^{k-1} (2^{k-j} - 1) \binom{k}{j}A_j, \\
     B_k \ &= \ \sum_{j=0}^{k-1} (2^{k-j} - 1) \binom{k}{j}B_j \ \ \ +  1. \label{HoggattB}
\end{align}
We note that Hoggatt's Equation (\ref{HoggattB}) follows from our Theorem \ref{iFi.theorem}, as we now demonstrate. If we take 
$n=2$ in Equation (\ref{iFi.equation}) 
from Theorem \ref{iFi.theorem}, we get
\begin{align*}
0^p\cdot F_2 + 1^p\cdot F_1 + 2^p\cdot F_0 \ &= \ A_p F_2 + B_p F_3 - \sum_{k=0}^p \binom{p}{k} B_k 2^{p-k}, 
\intertext{and if we simplify this and pull out the $k=p$ term from the sum, we get }
1 \ &= \ A_p + 2B_p - B_p - \sum_{k=0}^{p-1} \binom{p}{k} B_k 2^{p-k}.
\end{align*}
We now replace the $A_p$ in the above formula with the expression in Equation (\ref{DresdenA}) and simplify further to get
\[
B_p = 1 + \sum_{k=0}^{p-1} (2^{p-k} - 1) \binom{p}{k} B_k.
\]
If we replace $k$ with $j$, and $p$ with $k$, we obtain Hoggatt's formula in Equation (\ref{HoggattB}), as desired.

As alluded to earlier, Zeitlin \cite{Zeitlin} found  formulas for $A_k$ and $B_k$ in terms of factorials, Fibonacci numbers, and the Stirling numbers of the second kind. Here they are:
\begin{align}
     A_k \ &= \ \sum_{j=0}^{k} \, j! \, F_{j+1} 
                \begin{Bmatrix}k \\ j \end{Bmatrix}, \label{Zeitlin1}\\
     B_k \ &= \ \sum_{j=0}^{k} \, j! \, F_{j+2} 
                \begin{Bmatrix}k \\ j \end{Bmatrix}. \label{Zeitlin2}
\end{align}
Zeitlin also provided these lovely identities involving our numbers $A_k$ and $B_k$ along with the Stirling numbers of the first kind: 
\begin{align*}
     n!F_{n+1} \ &= \ \sum_{k=1}^{n} 
                \begin{bmatrix}n \\ k \end{bmatrix} A_k, \\
     n!F_{n+2} \ &= \ \sum_{k=1}^{n} 
                \begin{bmatrix}n \\ k \end{bmatrix} B_k. \\
     \end{align*}
     
Kinlaw, Morris, and Thiagarajan \cite{KMT} take a novel approach to defining $A_k$ and $B_k$. They use a  polylogarithm, along with the golden ratio $\phi$, to define $A_k$ and $B_k$ as
\[
A_k + B_k \phi \ = \ \sum_{n=1}^{\infty} \frac{n^k}{\phi^n}.
\]
From here, they find an expression for $A_k$ and $B_k$ in terms of the Fibonacci numbers and the Eulerian numbers  $\left\langle
\begin{matrix} k\\j \end{matrix}
\right\rangle$ as defined in \cite{GKP}. Their Proposition 4.1
states that 
\[
A_k = \sum_{j=0}^{k-1} 
\left\langle
\begin{matrix} k\\j \end{matrix}
\right\rangle
F_{k+j+1}, 
\qquad 
B_k = \sum_{j=0}^{k-1} 
\left\langle
\begin{matrix} k\\j \end{matrix}
\right\rangle
F_{k+j+2}. 
\]
Kinlaw, Morris, and Thiagarajan also show two more 
interesting properties. First, they demonstrate in their Proposition 4.4 that $A_k$ and $B_k$ are the coefficients in the exponential generating functions for $1/(e^x - e^{2x} + 1) - 1$ and 
$e^x/(e^x - e^{2x} + 1)$ respectively. Then, their Proposition 4.5 shows that their $A_k$ and $B_k$ numbers satisfty Zeitlin's formulas in Equations (\ref{Zeitlin1}) and (\ref{Zeitlin2}) with the Stirling numbers of the second kind, thus verifying that their $A_k$ and $B_k$ are the same as ours. 

Shannon and Ollerton \cite[Conjecture 5]{ShannonOllerton} posit the following relationship between our numbers $B_k$ for $k \geq 2$, and the Bernoulli numbers 
$\mathbb{B}_n^{+}$ where the $+$ means that the $n=1$ Bernoulli number is $+1/2$ instead of the typical $-1/2$. Here is their conjecture, written using our notation:
\[
B_k \ = \ k\Bigg(\frac{5}{2}B_{k-1} + (-1)^{k-1}\Bigg) \ - \  
\sum_{j=1}^{k-2} (-1)^{k+j} B_j \Bigg(\sum_{r=j}^{k} 
\binom{k}{r}\binom{r}{j} 
\mathbb{B}_{k-r}^{+}\Bigg).
\]
We have verified that this holds up to $k=1000$. 

Finally, Adegoke \cite[Theorem 1]{Adegoke} uses the derivatives of $1/(e^x - e^{2x} + 1)$ as seen in Kinlaw, Morris, and Thiagarajan  to establish the following nice property, which 
we have translated into our notation:
\[
B_k \ = \ (-1)^k \left[ 1 - \sum_{j=0}^{k-1} \binom{k}{j} (2^{k-j} +1) (-1)^j B_j \right].
\]
The attentive reader will notice that this is essentially an alternating version of Hoggatt's 
Equation (\ref{HoggattB}).

\section{Conclusion}

Why stop with the Fibonacci numbers? Both Ledin and Adegoke extended these ideas to the Lucas numbers, and those formulas are strikingly similar to the ones for the Fibonacci numbers. Zeitlin showed that similar formulas exist for arbitrary recurrence sequences, although he readily admits the difficulty of his technique when he writes \cite[p.~35]{Zeitlin}, ``A formidable obstacle in this procedure is the complex nature of the [coefficients], which involve multiple summations". 

Our insight is that we can start with simple formulas that involve only binomial coefficients, and then work our way towards the desired Brousseau sums.
As an encouragement for future work, 
we have this formula for the Pell numbers, defined 
as $P_0 = 0$, $P_1 = 1$, and $P_n = 2P_{n-1} + P_{n-2}$. Our formula is
\[
P_n - n^3 \ = \ 2 \sum_{i=1}^n \left(i^3 - 
\binom{3}{2} i^2  - \binom{3}{0}i^0\right)\cdot P_{n-i},
\]
and of course this is almost identical to our
nice little identity in Equation (\ref{nicelittle}) for the Fibonacci numbers, which is how we began this journey. We can even go one step further and consider the general linear recurrence defined as 
$R_0 = 0$, $R_1 = 1$, and $R_n = a R_{n-1} + b R_{n-2}$. It is not hard to show that 
\[
R_n - n^3 \ = \ \sum_{i=1}^n \left( (a+b-1)i^3 
- (b+1)\binom{3}{2} i^2 
+ (b-1)\binom{3}{1} i^1
- (b+1) \binom{3}{0}i^0\right)\cdot R_{n-i}.
\]
From this, one should be able to find a formula for $\sum_{i=1}^n i^3R_i$ and, more generally, $\sum_{i=1}^n i^p R_i$. We encourage the reader to continue on this path.

\end{document}